\documentclass[11pt]{article}
\usepackage{amssymb,amsmath,amsthm}
\usepackage{amscd}
\usepackage{mathrsfs}
\allowdisplaybreaks[1]
\newcommand{\N}{\mathbb{N}}
\newcommand{\R}{\mathbb{R}}

\newcommand{\dd}{\,{\rm d}}

\numberwithin{equation}{section}
\newtheorem{thm}{Theorem}[section]
\newtheorem{df}[thm]{Definition}
\newtheorem{prop}[thm]{Proposition}
\newtheorem{lem}[thm]{Lemma}
\newtheorem{rem}[thm]{Remark}
\newtheorem{cor}[thm]{Corollary}

\begin{document}

\title{Remark on the Helmholtz decomposition 
in domains with noncompact boundary}

\date{}

\author{
\null\\
Yasunori Maekawa\\    
Mathematical Institute, Tohoku University\\          
6-3 Aoba, Aramaki, Aoba, Sendai 980-8578, Japan\\
{\tt maekawa@math.tohoku.ac.jp }
\and
\\
Hideyuki Miura \\
Department of Mathematics,  Graduate School of Science,  Osaka University\\
1-1 Machikaneyama, Toyonaka, Osaka 560-0043, Japan\\
{\tt miura@math.sci.osaka-u.ac.jp}}

\maketitle

\begin{center}
{\bf Abstract}
\end{center}
\vspace{-2mm}

Let $\Omega$ be a domain with noncompact boundary. 
It is known that the Helmholtz decomposition is not 
always valid in  $L^p(\Omega)$ except for the energy space $L^2 (\Omega)$. 
In this paper we consider  a typical unbounded domain whose boundary is given as a Lipschitz graph, and show that 
the Helmholtz decomposition holds in certain anisotropic spaces which include some infinite energy vector fields.

\vspace{0.5cm}



% deduce, along with

\section{Introduction}\label{sec.intro}

%\subsection{Motivation and main results}

Let $\Omega \subset \R^{d+1}$ ($d \ge 1$) be 
a domain with Lipschitz boundary. 
The Helmholtz decomposition, the decomposition of a 
given vector field into a solenoidal field and a potential one is
the fundamental tool in the mathematical analysis of the incompressible flow. 
In the energy space $(L^2 (\Omega))^{d+1}$ 
this decomposition is easily derived for any domain $\Omega$ 
from the standard theory of the Hilbert space.  
On the other hand, if the space $(L^2 (\Omega))^{d+1}$  is replaced 
by other function spaces such as $(L^q (\Omega))^{d+1}$, then 
the verification of the Helmholtz decomposition requires 
detailed analysis in general. In the case when 
$\Omega$ is a bounded domain or an exterior domain 
with smooth boundaries, 
the validity of the decomposition in  $(L^q (\Omega))^{d+1}$, $1<q<\infty$, is shown by \cite{FM} and 
\cite{Miyakawa} respectively, and then 
their results are extended to these domains but with $C^1$-boundary 
by \cite{SS}. 
Moreover, for the bounded Lipschitz domains, 
the validity is proved around $3/2<q<3$ 
in \cite{FMM}, 
and for any $1<q<\infty$ by \cite{GZ} when the domain is convex.  
However, even if the boundary is smooth enough, the problem becomes subtle
when the boundary is noncompact. 
Although the decomposition is still valid for $1<q<\infty$
for some special cases, e.g., aperture domains \cite{Farwig},
layers \cite{Miyakawa94}, cylinders \cite{ST}, 
half spaces and their small perturbations \cite{SS}, 
it is known that the domain of simple form
\begin{equation}
\Omega = \{\tilde x = (x,x_{d+1})\in \R^d\times \R~|~x_{d+1}>\eta(x)\},
\label{domain}
\end{equation}
with a given function $\eta$ 
%satisfying  $\|\nabla_x\eta\|_{L^\infty (\R^d)}<\infty$
does not always admit the Helmholtz decomposition in $(L^q (\Omega))^{d+1}$ if $q\ne 2$, even if $\eta$ is smooth,
see \cite{Bogovskii} and \cite[III.1]{Galdi}.  
Hence it is an important question to ask which function space, 
other than 
 $(L^2 (\Omega))^{d+1}$, admits the Helmholtz decomposition.
In \cite{FKS,FKS2}, the authors considered $\tilde{L}^q(\Omega)$
defined by 
$$
\tilde{L}^q(\Omega)= 
\begin{cases}
&
L^2(\Omega) \cap L^q(\Omega)~~~~ \quad \ 2 \le q < \infty,
\\
&L^2(\Omega) + L^q(\Omega)~~~~~ \quad  1<q <2,
\end{cases}
$$
and showed that general domains with uniform $C^1$ boundaries
admit the Helmholtz decomposition in these spaces. 
In this paper, we will give an alternative approach
for this question in the domain 
of the form \eqref{domain}.
Before stating the result, it would be convenient to formulate 
our problem 
more systematically. Let $X(\Omega)$ be a Banach space of functions 
in $\Omega$ satisfying  $C_0^\infty (\Omega)\subset X(\Omega) 
\subset L^1_{loc}(\Omega)$. Set 
\begin{align}
X_\sigma (\Omega) = 
\overline{C_{0,\sigma}^\infty (\Omega )}^{\|\cdot\|_{X(\Omega)}},~~~~~~~~ 
X_G (\Omega ) = \{ \nabla p \in (X(\Omega ))^{d+1}~|~p\in  L^1_{loc}(\Omega) \}.
\end{align}
Here $C_{0,\sigma}^\infty (\Omega)$ is a set of all smooth, 
compactly-supported, and divergence-free vector fields in $\Omega$. 
For simplicity of the  notation we write $\|\cdot\|_{X(\Omega)}$
 for $\|\cdot \|_{(X(\Omega))^{d+1}}$.  
\begin{df}\label{def.helmholtz} We say that the space $(X(\Omega))^{d+1}$ 
admits the Helmholtz decomposition if each $f\in (X(\Omega))^{d+1}$ 
has a unique decomposition  $f=u+\nabla p$, $u\in X_\sigma (\Omega)$, 
$\nabla p\in X_G(\Omega)$, satisfying
\begin{align}
\| u\|_{X(\Omega)} + \|\nabla p \|_{X(\Omega)}\leq C \| f\|_{X(\Omega)}.\label{eq.def.helmholtz}
\end{align}
Here $C$ is a positive constant independent of $f$. 
\end{df}

In order to consider the domain $\Omega$ of the form \eqref{domain} we 
define the standard isomorphism $\Phi:\Omega\in \tilde x  \mapsto 
\tilde y = \Phi(\tilde x) \in \R^{d+1}_+$  by
\begin{align}\label{def.diffeo}
\Phi_j (\tilde x) =
\begin{cases}
& x_j~~~~~~~~~~~~~~~~~~~~{\rm  if}~~ 1\leq j\leq d,\\
& x_{d+1}-\eta (x)~~~~~~~~{\rm if}~~j=d+1.
\end{cases}
\end{align} 
 Let $1<q,r<\infty$ and let $Y^{q,r}(\Omega)$ be the Banach space defined by 
\begin{align} 
Y^{q,r}(\Omega) = \{ f\in L^1_{loc}(\Omega)~|~\| f\|_{Y^{q,r}(\Omega)} = \| f\circ \Phi^{-1} \|_{L^q_t  (\R_+; L^r_y (\R^d) )} <\infty \}
\end{align}
with the norm $\|\cdot \|_{Y^{q,r}(\Omega)}$. Here we have used the notation $\tilde y = (y,t)\in \R^d\times \R_+$. Our main result reads as follows:
\begin{thm}
\label{thm.helmholtz.intro}
Let $\Omega$ be a domain of the form \eqref{domain} with uniform Lipschitz
boundary. Then the space $\big (Y^{q,2}(\Omega) \big )^{d+1}$ 
admits the Helmholtz decomposition for all $1<q<\infty$. 
Moreover, the constant $C$ in \eqref{eq.def.helmholtz} depends only on $d$, 
$q$, and $\|\nabla_x \eta\|_{L^\infty (\R^d)}$. 
\end{thm}

\begin{rem}{\rm By the definition it is easy to see 
that the dual space of $Y^{q,r}(\Omega)$, 
denoted by $Y^{q,r}(\Omega)^*$, 
is the space 
$Y^{q',r'}(\Omega)$ with $q'=q/(q-1)$ and $r'=r/(r-1)$. 
The space $(Y^{q,q}(\Omega))^{d+1}$ coincides with $(L^q(\Omega))^{d+1}$. However, due to the well-known counterexample 
of the weak Neumann problem in the exterior of the cone-like domain \cite{Bogovskii},
one cannot expect the validity of the Helmholtz decomposition
in the usual $L^q$ space. 
On the other hand, compared with the result in \cite{FKS}, 
our function space $Y^{q,2}(\Omega)$ includes a class of functions 
decaying slowly in the $x_{d+1}$ direction  since Theorem \ref{thm.helmholtz.intro} allows any $q\in (2,\infty)$.

%It should be emphasized that no smallness  condition is required for $\|\nabla_%x \eta\|_{L^\infty (\R^d)}$ in Theorem \ref{thm.helmholtz.intro},
%and any uniform Lipshitz domains 
%of the form \ref{domain} are allowed.
}
\end{rem}
As is well-known, the verification of the Helmholtz decomposition is reduced to the unique solvability of the weak Neumann problem in 
$(Y^{q,2}(\Omega))^{d+1}$ :
\begin{align}\label{weak.neumann.helmholtz}
\langle \nabla p, \nabla \varphi\rangle _{L^2(\Omega)} = \langle f, \nabla \varphi\rangle _{L^2(\Omega)} ~~~~~~~~ {\rm for~all}~~\varphi\in \{ f\in L^1_{loc}(\Omega)~|~\nabla f\in \big (Y^{q',r'}(\Omega )\big )^{d+1}\},
\end{align}
which is a weak formulation of 
\begin{align}\label{neumann.helmholtz}
\Delta p = \nabla\cdot f~~~~{\rm in }~~\Omega, ~~~~~~~~~~~~~~n\cdot \nabla p = n\cdot f~~~~{\rm on}~~\partial\Omega.
\end{align}
Here $n$ stands for the exterior unit normal to $\partial\Omega$. 
Through the isomorphism $\Phi$ defined by \eqref{def.diffeo}, 
the problem \eqref{weak.neumann.helmholtz} or 
\eqref{neumann.helmholtz} is transformed to the Neumann problem for 
an elliptic partial differential equation of the divergence form 
in $\R^{d+1}_+$. 
%with a matrix $A$ which is real and positive symmetric. 
Then our task is to look for the solution of the transformed problem in 
$L^q(\R_+; L^2(\R^d))$. 
For the purpose, we will make use of an approach 
proposed in the companion work \cite{MM},
where we gave a solution formula
for the boundary value problem to divergence form elliptic equations 
in $\R^{d+1}_+$ in terms of the Poisson semigroups; see Theorem \ref{thm.factorization} for details. This solution formula combined with the semigroup theory 
 yields a sufficient condition for function spaces
to ensure the solvability of the Neumann problem, and it will be verified that $Y^{q,2}$ satisfies this condition.

Before concluding the introduction, we would like to point out 
the difference of our approach from previous works on the Neumann problem 
in the domain with regular (e.g. uniformly $C^1$) boundary. 
In \cite{SS, FKS, GHHS}, they employed localization
procedure for the Neumann problem 
to reduce the problem to a countable number of the Neumann problems in 
$\R^{d+1}_+$ or $\R^{d+1}$. Then thanks to the regularity assumption
of the boundary, each problem can be dealt with as 
a small perturbation of the Neumann problem for the Poisson equation.
On the other hand, our approach is not relied on this perturbation technique,
and therefore one can handle even \textit{large} perturbation
of the Neumann problem 
for the Poisson equation 
in $\R^{d+1}_+$ 
with respect to the Lipschitz norm of  $\eta$.
Instead, we need to use the $L^2$ space in the $x$ direction. 

In the next section we will recall the solution formula for the Neumann problem,
 which is 
based on the factorization of the elliptic operators 
in \cite{MM}.
Then we will prove the main theorem in Section 3.

\section{Solution formula for the Neumann problem in $\R^{d+1}_+$} 

Consider the second order elliptic operator of divergence form in $\R^{d+1}= \{ (x,t) \in \R^d \times \R \}$, 
\begin{equation}
\mathcal{A} = -\nabla \cdot A \nabla, ~~~~~~~~~~~~~~A = A (x) = \big (a_{i,j} (x) \big ) _{1\leq i,j\leq d+1}~.\label{def.calA}
\end{equation}
Here $d\in \N$, $\nabla = (\nabla_x,\partial_t)^\top$ with $\nabla_x=(\partial_1,\cdots,\partial_d)^\top$, and each $a_{i,j}$ is always assumed to be 
\textit{$t$-independent}. We further assume that $A$ is a \textit{real symmetric} matrix and each component $a_{i,j}$ 
is a measurable function satisfying the
uniformly elliptic condition
\begin{align}
\langle A (x) \eta, \eta \rangle \geq \nu_1 |\eta|^2, ~~~~~~~~~~ | \langle A (x) \eta,\zeta\rangle | \leq \nu_2 |\eta| |\zeta|\label{ellipticity}
\end{align}
for all $\eta,\zeta\in \R^{d+1}$ and for some constants $\nu_1,\nu_2$ with $0<\nu_1\leq \nu_2 <\infty$. Here $\langle \cdot,\cdot \rangle$ denotes the inner product of $\R^{d+1}$, i.e., $\langle \eta, \zeta\rangle = \sum_{j=1}^{d+1}\eta_j \zeta_j$ for $\eta,\zeta\in \R^{d+1}$. For later use we set $b=a_{d+1, d+1}$, which satisfies $\nu_1\leq b\leq \nu_2$ due to \eqref{ellipticity}. We also denote by ${\bf a}$ the vector ${\bf a} (x)  = ( a_{1, d+1} (x) ,\cdots , a_{d, d+1}(x))^\top$. 

We denote by $D_H(T)$ the domain of a linear operator $T$ in a Banach space $H$. 
Under the condition \eqref{ellipticity} the standard theory of sesquilinear 
forms gives a realization of $\mathcal{A}$ in $L^2 (\R^{d+1})$, 
denoted again by $\mathcal{A}$, such as \begin{align}
D_{L^2}(\mathcal{A}) & = \big \{ w \in H^1 (\R^{d+1})~|~{\rm there ~is}~ F \in L^2 (\R^{d+1})~{\rm such ~that}~\nonumber \\
& ~~~~~ \langle A\nabla w, \nabla v\rangle _{L^2(\R^{d+1})} = \langle F, v\rangle _{L^2 (\R^{d+1})}~{\rm for~all}~v\in H^1 (\R^{d+1})\big \},
\end{align}
and $\mathcal{A} w  = F$ for $w\in D_{L^2}(\mathcal{A})$. Here $H^1(\R^{d+1})$ is the usual Sobolev space and $\langle w,v\rangle_{L^2(\R^{d+1})}= \int_{\R^{d+1}} w (x,t) v (x,t) \dd x \dd t$. 
\begin{df}\label{def.intro} {\rm (i)} For a given $h\in \mathcal{S}' (\R^d)$ we denote by $M_h : \mathcal{S}(\R^d)\rightarrow \mathcal{S}'(\R^d)$ the multiplication $M_h u = h u$.

\noindent {\rm (ii)} We denote by $E_{\mathcal{A}}: {H}^{1/2} (\R^d)\rightarrow \dot{H}^1 (\R^{d+1}_+)$ the $\mathcal{A}$-extension operator, i.e., $w=E_{\mathcal{A}} \varphi$ is the solution to the Dirichlet problem
\begin{equation}\label{eq.dirichlet0}
\begin{cases}
& \mathcal{A} w  = 0~~~~~~~{\rm in}~~~\R^{d+1}_+,\\
& \hspace{0.3cm} w  = \varphi~~~~~~~{\rm on}~~\partial\R^{d+1}_+=\R^d.
\end{cases}
\end{equation}
The one parameter family of linear operators $\{E_{\mathcal{A}} (t)\}_{t\geq 0}$, defined by $E_{\mathcal{A}} (t) \varphi = (E_{\mathcal{A}} \varphi )(\cdot,t)$ for $\varphi \in H^{1/2}(\R^d)$, is called the Poisson semigroup associated with $\mathcal{A}$.
 
\noindent {\rm (iii)} We denote by $\Lambda_{\mathcal{A}}: {H}^{1/2} (\R^d)\rightarrow \dot{H}^{-1/2} (\R^d)=\big (\dot{H}^{1/2} (\R^d) \big )^*$ the Dirichlet-Neumann map associated with $\mathcal{A}$, which is defined through the sesquilinear form
\begin{equation}
\langle  \Lambda_{\mathcal{A}} \varphi, g \rangle_{\dot{H}^{-\frac12},\dot{H}^{\frac12}}  =  \langle  A\nabla E_{\mathcal{A}}\varphi, \nabla E_{\mathcal{A}} g \rangle_{L^2(\R^{d+1}_+)}, ~~~~~~~~~~\varphi,g \in H^{1/2}(\R^d).\label{def.Lambda}
\end{equation}
Here $\langle \cdot,\cdot\rangle _{\dot{H}^{-1/2},\dot{H}^{1/2}}$ denotes the duality coupling of  $\dot{H}^{-1/2}(\R^d)$ and $\dot{H}^{1/2}(\R^d)$.

\end{df}
\begin{rem}{\rm (i) As usual, Eq. \eqref{eq.dirichlet0} is considered in a weak sense; cf. \cite[Section 2.1]{MM}. The proof of the existence of the extension operator $E_{\mathcal{A}}$ is 
classical, and indeed it is a consequence of the Riesz representation 
theorem together with the harmonic extension of the function in 
$H^{1/2}(\R^d)$. As is shown in \cite[Proposition 2.4]{MM},  $\{E_{\mathcal{A}} (t)\}_{t\geq 0}$ is a strongly continuous and analytic semigroup in $H^{1/2}(\R^d)$. We denote its generator by $-\mathcal{P}_{\mathcal{A}}$, and $\mathcal{P}_{\mathcal{A}}$ is called a {\it Poisson operator} associated with $\mathcal{A}$. (ii) Since $A$ is Hermite and satisfies the uniformly elliptic condition
\eqref{ellipticity}, 
the theory of the sesquilinear forms  \cite[Chapter VI. 
\S 2]{Kato}
shows that  $\Lambda_{\mathcal{A}}$ is extended as 
a self-adjoint operator in $L^2(\R^d)$.
}
\end{rem}
 
The following result plays a fundamental role in the derivation of the solution formula for the Neumann problem.
\begin{thm}[{\cite[Theorem 1.3, Theorem 4.2]{MM}}]\label{thm.factorization} \, Let $\mathcal{A}$ be the elliptic operator defined in \eqref{def.calA}
with a real symmetric matrix  $A$ satisfying \eqref{ellipticity}.
Then $D_{L^2}(\Lambda_{\mathcal{A}})=H^1 (\R^d)$ with equivalent norms and the operator $-{\bf P}_{\mathcal{A}}$  defined by  
\begin{align}
D_{L^2} ({\bf P}_{\mathcal{A}}) = H^1 (\R^d), ~~~~~~~~~ -{\bf P}_{\mathcal{A}} \varphi =  -M_{1/b} \Lambda_{\mathcal{A}} \varphi - M_{{\bf a}/b} \cdot \nabla _x \varphi, 
\label{def.P}
\end{align}
generates a strongly continuous and bounded analytic semigroup  in $L^2 (\R^d)$. Moreover, the realization 
$\mathcal{A}'$ in $L^2 (\R^d)$ and the realization $\mathcal{A}$ in $L^2 (\R^{d+1})$ are respectively factorized as 
\begin{align}
\mathcal{A}' & = M_b \mathcal{Q}_{\mathcal{A}} {\bf P}_{\mathcal{A}}, ~~~~~~~~~ \mathcal{Q}_{\mathcal{A}} = M_{1/b}  ( M_b {\bf P}_{\mathcal{A}} )^*,\label{eq.factorization'}\\
\mathcal{A} & = - M_b (\partial_t -  \mathcal{Q}_{\mathcal{A}} ) (\partial_t + {\bf P}_{\mathcal{A}} ).\label{eq.factorization}
\end{align}
Here $(M_b {\bf P}_{\mathcal{A}})^*$ is the adjoint of $M_b {\bf P}_{\mathcal{A}}$ in $L^2 (\R^d)$.  

\end{thm}

\begin{rem}{\rm  
The operator ${\bf P}_{\mathcal{A}}$ is nothing but the Poisson operator $\mathcal{P}_{\mathcal{A}}$ associated with $\mathcal{A}$. That is, $-{\bf P}_{\mathcal{A}} \varphi = -\mathcal{P}_{\mathcal{A}} \varphi := \lim_{t\downarrow 0} t^{-1} \big ( E_{\mathcal{A}} (t) \varphi   - \varphi \big )$ in $L^2 (\R^d)$ for $\varphi \in H^1 (\R^d)$.
}
\end{rem}

Now we  consider the inhomogenuous Neumann problem
\begin{equation}\label{eq.neumann}
\begin{cases}
& ~~~~~~~~~~~\hspace{0.3cm}~ \mathcal{A} w  = F~~~~~~~{\rm in}~~~\R^{d+1}_+,\\
& -\langle {\bf e}_{d+1}, A \nabla w \rangle  = g~~~~~~~{\rm on}~~\partial\R^{d+1}_+.
\end{cases}
\end{equation}
By a direct application of the factorization
\eqref{eq.factorization}, 
one can easily derive the formal representation 
of the solution to the Neumann problem \eqref{eq.neumann} 
as follows.
\begin{thm}[{\cite[Theorem 5.1]{MM}}]\label{thm.representation.formula} Assume that $F,\partial_t F \in \dot{H}^{-1} (\R^{d+1}_+)$, $g\in H^{1/2}(\R^d)$. Assume further that $h=g + M_b \int_0^\infty e^{-s\mathcal{Q}_{\mathcal{A}}} M_{1/b} F (s) \dd s$ belongs to the range of $\Lambda_{\mathcal{A}}$ in $L^2 (\R^d)$. Let $w \in C([0,\infty); L^2 (\R^d))\cap \dot{H}^1(\R^{d+1}_+)$ be a weak solution to \eqref{eq.neumann}. Then $w$ is a mild solution, i.e, 
\begin{align}
& w (t) =  e^{-t \mathcal{P}_{\mathcal{A}}} \Lambda_{\mathcal{A}}^{-1} \big ( g + M_b   \int_0^\infty   e^{-s \mathcal{Q}_{\mathcal{A}}} M_{1/b} F  (s) \dd s \big )   + \int_0^t e^{-(t-s) \mathcal{P}_{\mathcal{A}}} \int_s^\infty  e^{-(\tau-s) \mathcal{Q}_{\mathcal{A}}} M_{1/b} F  (\tau)  \dd \tau \dd s.
\label{eq.formula1}
\end{align}
\end{thm}
%This formula makes sense at least if the data $f$, $g$ are smooth
%and if
%\begin{equation}
%h:
%=g + M_b \int_0^\infty e^{-s \mathcal{Q}_{\mathcal{A}}} M_{1/b} f (s) \dd s
%\end{equation}
%belongs to the range of $\Lambda_A$
We note that $e^{-t\mathcal{Q}_{\mathcal{A}}}$ is related with $e^{-t\mathcal{P}_{\mathcal{A}}}$ through the formula 
\begin{align}
e^{-t\mathcal{Q}_{\mathcal{A}}} = M_{1/b} \big ( e^{-t\mathcal{P}_{\mathcal{A}}} \big )^* M_b.
\end{align}
Then the representation \eqref{eq.formula1} reduces the inhomogeneous problem \eqref{eq.neumann} to the analysis of the semigroup $\{e^{-t\mathcal{P}_{\mathcal{A}}}\}_{t\geq 0}$ and the operator $\Lambda_{\mathcal{A}}$.

 %Theorem \ref{thm.factorization} is now an important tool to prove the following result.

%\begin{thm}\label{thm.kato.problem} Assume that the semigroup $\{e^{t\mathcal{P}_{\mathcal{A}}}\}_{t\geq 0}$ acting on $H^{1/2}(\R^d)$ is extended to the strongly continuous semigroup acting on $L^2 (\R^d)$. Then   $D_{L^2}(\mathcal{P}_{\mathcal{A}})$ is continuously embedded in $H^1 (\R^d)$. Moreover, if the extended semigroup is bounded, i.e., $\displaystyle \sup_{t>0} \| e^{t\mathcal{P}_{\mathcal{A}}} \|_{L^2\rightarrow L^2}<\infty$, then we have 
%\begin{align}\label{eq.thm.kato.problem}
%\| \nabla_x u \|_{L^2(\R^d)} \leq C \| \mathcal{P}_{\mathcal{A}} u\|_{L^2(\R^d)}~~~~~~~{\rm  for~ all~~} u\in D_{L^2}(\mathcal{P}_{\mathcal{A}}).
%\end{align}
  
%\end{thm}

\section{Helmholtz decomposition in $(Y^{q,r}(\Omega))^{d+1}$}\label{sec.proof}

As stated in the introduction, the Helmholtz decomposition for a given vector field is reduced to the Neumann problem \eqref{neumann.helmholtz}. 
Let $\Phi:\Omega\rightarrow \R^{d+1}_+$ be the isomorphism defined by \eqref{def.diffeo}. By taking the push-forward 
\begin{align}
w = p \circ \Phi^{-1},~~~~~~~~~F=(F',F_{d+1})=f\circ \Phi^{-1}, \label{def.pushforward}
\end{align}
the problem \eqref{neumann.helmholtz} is transformed to the Neumann problem in $\R^{d+1}_+$:
\begin{equation}\label{eq.helmholtz}
\begin{cases}
\quad\quad\quad\quad \mathcal{A} w  
& =  
-\nabla_x\cdot F' -  \partial_t(F_{d+1}+ M_{\bf a} \cdot F')
~~~~~{\rm in}~~\R^{d+1}_+,\\
-\langle {\bf e}_{d+1}, A\nabla w\rangle  
& = 
- ( F_{d+1} + M_{{\bf a}} \cdot F')
~~~~~~~~~~~~~~~~~~~~{\rm on}~~\partial\R^{d+1}_+.
\end{cases}
\end{equation}
Here the matrix $A$ in this case is real symmetric and positive definite with
 ${\bf a} = -\nabla_x \eta$, $b = 1+|\nabla_x \eta |^2$, and $A'=(a_{i,j})_{1\leq i,j\leq d}=I'$ (the identity matrix). Let $1<q,r<\infty$ and set 
\begin{align}
Z^{q,r}(\R^{d+1}_+) :&=\dot{W}^{1,q}(\R_+; L^{r}(\R^d))
\cap L^{q}(\R_+; \dot{W}^{1,r}(\R^d))
\notag
\\
&= \{ \phi \in L^1_{loc}(\R^{d+1}_+)~|~\partial_i \phi \in L^{q}(\R_+; L^{r}(\R^d) ) ~~~ 1\leq i\leq d+1 \}.
\end{align} 
Let $F\in  \big (L^{q}(\R_+; L^{r}(\R^d) )\big )^{d+1}$. The weak formulation of \eqref{eq.helmholtz} is then to look for $w\in Z^{q,r}(\R^{d+1}_+)$ such that 
\begin{align}
\langle A\nabla w, \nabla \phi\rangle _{L^2(\R^{d+1}_+)} = \langle F', \nabla_x \phi + M_{\bf a} \partial_t\phi \rangle _{L^2(\R^{d+1}_+)} + \langle F_{d+1}, \partial_t \phi\rangle _{L^2(\R^{d+1}_+)}\label{eq.helmholtz.weak}
\end{align}
for all $\phi \in  Z^{q',r'}(\R^{d+1}_+)$. 

In the following paragraphs we abbreviate $\mathcal{P}_{\mathcal{A}}$ ($\mathcal{Q}_\mathcal{A}$, $\Lambda_\mathcal{A}$) to $\mathcal{P}$
($\mathcal{Q}$ and $\Lambda$ as well) for simplicity of the notation. The most important step in the analysis of \eqref{eq.helmholtz} is to derive the estimate corresponding 
with \eqref{eq.def.helmholtz}, which is closely related to the spectral properties of $\mathcal{P}$ and $\Lambda$. 
To make the essence of our arguments clear, we will give  in Section \ref{subsec.helmholtz.L^r} natural sufficient conditions for the Helmholtz decomposition to hold in $(Y^{q,r}(\Omega))^{d+1}$ in terms of the properties of  $\mathcal{P}$ and $\Lambda$ in $L^r(\R^d)$. Roughly speaking, the following three conditions are required: Let $1<m<\infty$ and set $m'=m/(m-1)$. (i) boundedness of the semigroups $\{e^{-t\mathcal{P}}\}_{t\geq 0}$,  $\{e^{-t\Lambda }\}_{t\geq 0}$ in $L^r(\R^d)$, $r=m,m'$, (ii) coercive estimates for $\mathcal{P}$ and $\Lambda$ in $L^r (\R^d)$, $r=m,m'$, (iii) maximal regularity estimates for $\{e^{-t\mathcal{P}}\}_{t\geq 0}$ in $L^r (\R_+; L^r (\R^d))$, $r=m,m'$. As long as $\nabla_x\eta$ is uniformly bounded, these conditions are shown to hold at least for  $m=2$, which leads to Theorem \ref{thm.helmholtz.intro}; see Section \ref{subsec.helmholtz.L^2}.

\subsection{Sufficient condition for solvability of 
the Neumann problem}
\label{subsec.helmholtz.L^r}

In this section we investigate the relation between the boundedness of the Helmholtz decomposition and the spectral properties of  $\mathcal{P}$, $\Lambda$.  
%For a while we will assume that $\nabla_x \eta\in {\rm Lip} (\R^d)$. Then we can apply the results of \cite[Theorems 1.2, 1.3]{MM2}, and $D_{L^2}(\mathcal{P})= D_{L^2}(\Lambda)= D_{L^2}(\mathcal{Q}) = H^1(\R^d)$ holds with equivalent norms.
Let $1<m<\infty$ and $m'=m/(m-1)$. Let us recall that both $\mathcal{P}$  and  $\Lambda$ generate strongly continuous and bounded analytic seimgroups in $L^2 (\R^d)$. To develop our argument within the general  $L^r$ framework we first assume that 

\vspace{0.3cm}

\noindent {\bf (i)} The restrictions of  $\{e^{-t\mathcal{P}}\}_{t\geq 0}$ and $\{e^{-t\Lambda }\}_{t\geq 0}$ on $L^2(\R^d)\cap L^r(\R^d)$ are extended as  strongly continuous and bounded semigroups in $L^r (\R^d)$ with $r=m,m'$, i.e.,
\begin{equation}
\| e^{-t\mathcal{P}} \varphi \|_{L^r(\R^d)}  + \| e^{-t\Lambda} \varphi  \|_{L^r(\R^d)} \leq C \| \varphi \|_{L^r(\R^d)},~~~~~~~~~~~ t>0,~~~ \varphi \in L^r(\R^d).\label{condition2'}
\end{equation}

In fact,  the statement in {\bf (i)} is always verified at least for $\{e^{-t\Lambda}\}_{t\geq 0}$. While the behavior of  the Poisson semigroup in $L^r (\R^d)$ seems to be more difficult to analyze, and the estimate \eqref{condition2'}  for   $\{e^{-t\mathcal{P}}\}_{t\geq 0}$ can be obtained  at least for the case $r\in [2,\infty)$. We will sketch their proofs in the appendix for reader's convenience. In addition to {\bf (i)} we  assume in this section the following estimates {\bf (ii)} - {\bf (iii)}:

\vspace{0.3cm}  

\noindent {\bf (ii)} Coercive estimates: Let $r=m,~m'$. Then $D_{L^r}(\mathcal{P}) \cup D_{L^r}(\Lambda) \subset W^{1,r} (\R^d)$ and 
\begin{align}
&\|\nabla_x \varphi\|_{L^r(\R^d)} \leq C \| \mathcal{P} \varphi \|_{L^r(\R^d)},~~~~~~\varphi \in D_{L^r}(\mathcal{P}),\label{condition1-1}\\
&\| \nabla_x \varphi \|_{L^r(\R^d)} \leq C \|  \Lambda \varphi  \|_{L^r(\R^d)},~~~~~~\varphi \in D_{L^r}(\Lambda).\label{condition1-2}
\end{align}

\noindent{\bf (iii)} Maximal regularity: Let $r=m,~m'$.  The function $\Psi_\mathcal{P} [\phi](t) = \int_0^t e^{-(t-s)\mathcal{P}} \phi (s) \dd s$ satisfies 
\begin{align}
\| \mathcal{P} \Psi_{\mathcal{P}} [\phi ] \|_{L^r(\R_+; L^r (\R^d))} \leq C \| \phi \|_{L^r(\R_+; L^r(\R^d))}, ~~~~ \phi \in L^r (\R^{d+1}_+).\label{condition3}
\end{align}

\begin{rem}\label{rem.thm.helmholtz.1.1}{\rm As is well-known,  \eqref{condition2'} and \eqref{condition3} imply the analyticity of $\{e^{-t\mathcal{P}}\}_{t\geq 0}$ in $L^r(\R^d)$:
\begin{align}
\| t\mathcal{P} e^{-t\mathcal{P}} \varphi \|_{L^r(\R^d)} \leq C \| \varphi \|_{L^r(\R^d)},~~~~~~~~~t>0,~~ ~\varphi \in L^r (\R^d). \label{condition2}
\end{align}

}
\end{rem}

\begin{rem}\label{rem.thm.helmholtz.1.2}
{\rm Set $e^{-t\mathcal{P}}=0$ for $t<0$ and define the operator $\tilde \Psi_{\mathcal{P}}$ by $\tilde \Psi_\mathcal{P} [\phi] (t) = \int_\R e^{-(t-s)\mathcal{P}} \phi(s) \dd s$. Then \eqref{condition3} implies the estimate
\begin{align}
\| \mathcal{P} \tilde \Psi_{\mathcal{P}} [\phi] \|_{L^r (\R; L^r(\R^d))} \leq C \| \phi\|_{L^r(\R; L^r(\R^d))}.\label{condition3'}
\end{align}
From \eqref{condition2} and \eqref{condition3'} the theory of singular integral operators \cite{BCP} implies that
\begin{align}
\| \mathcal{P} \tilde \Psi_{\mathcal{P}} [\phi] \|_{L^q(\R; L^r (\R^d))} \leq C_q \| \phi \|_{L^q(\R; L^r(\R^d))},~~~~~~~~~~1<q<\infty.\label{condition4}
\end{align}

}
\end{rem}

\begin{rem}
\label{rem.inverse}
{\rm By the assumption {\bf (i)},  $\mathcal{P}$ and $\Lambda $ are sectorial in $L^r(\R^d)$ in the sense of \cite[Chapter 2]{Haase}. Thus 
the decompositions $L^r(\R^d)= 
{\rm Ker} (\mathcal{P})
\oplus 
\overline{{\rm Ran}(\mathcal{P})}$ and $L^r(\R^d) = {\rm Ker} (\Lambda)
\oplus 
\overline{{\rm Ran}(\Lambda)} $ hold;
see \cite[Proposition 2.2.1]{Haase} for the proof.
Moreover, since the operators are injective by {\bf (i)}, 
we see $L^r(\R^d)=\overline{{\rm Ran}(\mathcal{P})}=
\overline{{\rm Ran}(\Lambda)}$ and  the inverse operators $\mathcal{P}^{-1}$, $\Lambda^{-1}$
 can be extended to bounded operators from 
$L^r (\R^d)$ to the homogeneous Sobolev space $\dot{W}^{1,r}(\R^d)$.} 
\end{rem}

\begin{rem}{\rm In Section \ref{subsec.helmholtz.L^2} we will see that {\bf (ii)} and  {\bf (iii)} are always satisfied at least for $r=2$. 

}
\end{rem}

\begin{prop}
\label{prop.helmholtz}
Assume that  {\rm {\bf (i)} - {\bf (iii)}} hold. Let $F\in \big (C_0^\infty (\R^{d+1}_+) \big )^{d+1}$. Then there exists a unique weak 
solution $w\in \dot{H}^1 (\R^{d+1}_+)$ to \eqref{eq.helmholtz} satisfying 
\begin{align}
\| \nabla w \|_{L^q (\R_+; L^r (\R^d))} \leq C \| F \|_{L^q (\R_+; L^r (\R^d))},~~~~~~~~~1<q<\infty,~~\min\{m,m'\}\leq r\leq \max\{m,m'\}.\label{est.thm.helmholtz.1}
\end{align} 
Here $C$ depends only on $m$, $q$, $d$, $\|\nabla_x \eta\|_{L^\infty(\R^d)}$, and the constants in the estimates of  {\rm {\bf (i)} - {\bf (iii)}}.  
\end{prop}

Note that it suffices to show \eqref{est.thm.helmholtz.1} for $r=m,m'$ by the interpolation.  We start from the following lemma.
\begin{lem}\label{lem.helmholtz.1} 
Assume that  {\rm {\bf (i)} - {\bf (iii)}} hold. Set $e^{-t\mathcal{Q}} =0$ for $t<0$. Then $\tilde \Psi_{\mathcal{Q}} [\phi] (t)= \int_\R e^{-(t-s)\mathcal{Q}} \phi (s) \dd s$ satisfies 
\begin{align}
\| \mathcal{Q} \tilde \Psi_{\mathcal{Q}} [\phi] \|_{L^q (\R; L^r (\R^d))} \leq C_q \| \phi\|_{L^q (\R; L^r (\R^d))},~~~~~~~~~~1<q<\infty,~~r=m,~m'.\label{est.lem.helmholtz.1}
\end{align}

\end{lem}

\noindent {\it Proof.} We appeal to the duality argument. Since  $e^{-t\mathcal{Q}} = M_{1/b} (e^{-t\mathcal{P}})^* M_{b}$, we have for any $\psi \in C_0^\infty (\R^{d+1})$,
\begin{align}
\langle \mathcal{Q} \tilde \Psi_{\mathcal{Q}} [\phi] , \psi \rangle _{L^2 (\R^{d+1})} & = \int_\R \int_\R \langle  \mathcal{Q} e^{-(t-s)\mathcal{Q}} \phi (s),  \psi (t) \rangle _{L^2(\R^d)} \dd s \dd t \nonumber \\
& =\langle M_{b} \tilde \phi,  \mathcal{P}\tilde \Psi_{\mathcal{P}} [M_{1/b}\tilde \psi] \rangle_{L^2(\R^{d+1})}. \label{proof.lem.helmholtz.1.1}
\end{align}
Here $\tilde \phi (t) = \phi (-t)$ and $\tilde \psi  (t)  = \psi (-t)$. Then \eqref{est.lem.helmholtz.1} 
follows from \eqref{condition4} and the duality. 
The proof is complete.
\hfill
$\square$

\begin{lem}\label{lem.thm.helmholtz.2} Assume that {\bf (i)} - {\bf (iii)} hold. Then $D_{L^r}(\Lambda) \subset D_{L^r} (\mathcal{P})$, and  $M_{1/b}\nabla_x\cdot F'$ and $M_b \mathcal{Q} e^{-t\mathcal{Q}} F_{d+1}$ respectively belong to ${\rm Ran} (\mathcal{Q})$ and ${\rm Ran} (\Lambda)$ in $L^r(\R^d)$ for any  $F=(F',F_{d+1})\in (C_0^\infty (\R^d))^d\times C_0^\infty (\R^d)$ and $t>0$. Moreover,  it follows that
\begin{align}
&\| \mathcal{Q}^{-1}  M_{1/b}\nabla_x\cdot F'\|_{L^r (\R^d)}  + \| M_b  \mathcal{P} \Lambda^{-1}  F_{d+1} \|_{L^r (\R^d)}  \nonumber  \\
&~~~~~~~~~~~~~~~~~~~~~~~~~~~~~~~~~~~~~~~~+ \| \Lambda^{-1} M_b \mathcal{Q} F_{d+1} \|_{L^r (\R^d)} \leq C \| F\|_{L^r(\R^d)},~~~~~~~ r=m,~m'.
\end{align}
Here 
\begin{align*}
M_b  \mathcal{P} \Lambda^{-1}  F_{d+1} &  : = \displaystyle \lim_{\lambda\downarrow 0} M_b  \mathcal{P} (\Lambda+\lambda)^{-1}  F_{d+1}~~~~~~{\rm  in}~~ L^r(\R^d),\\
\Lambda^{-1}M_b \mathcal{Q}  F_{d+1} &  : = \displaystyle \lim_{t\rightarrow 0} \Lambda^{-1} M_b \mathcal{Q} e^{-t\mathcal{Q}} F_{d+1}~~~~~~{\rm  in}~~ L^r(\R^d).
\end{align*}

\end{lem}

\noindent {\it Proof.} We first prove that $M_{1/b}\nabla_x\cdot F'$ belongs to ${\rm Ran}(\mathcal{Q)}$ in $L^r(\R^d)$ and $\mathcal{Q}^{-1}  M_{1/b}\nabla_x\cdot$ is extended as a bounded operator from $(L^r(\R^d))^d$ to $L^r(\R^d)$. To this end we take any $\lambda>0$ and $\varphi\in C_0^\infty (\R^d)$ and then use the relation $e^{-t\mathcal{Q}}=M_{1/b} (e^{{-t\mathcal{P}}})^*M_{b}$  to derive 
\begin{align*}
\langle (\mathcal{Q} + \lambda)^{-1}M_{1/b} \nabla_x\cdot F', \varphi \rangle _{L^2(\R^d)} &= - \langle F', \nabla_x (\mathcal{P}+\lambda)^{-1} M_{1/b} \varphi \rangle _{L^2(\R^d)}\\
& = - \langle F', \nabla_x \mathcal{P}^{-1} \mathcal{P} (\mathcal{P}+\lambda)^{-1} M_{1/b} \varphi \rangle _{L^2(\R^d)}.
\end{align*}
Since $-\nabla_x\mathcal{P}^{-1}: {\rm Ran} (\mathcal{P})\rightarrow L^{r'}(\R^d)$ is extended  as a bounded operator $K$ from $L^{r'}(\R^d)$ to $(L^{r'}(\R^d))^d$ by the assumptions (cf. Remark \ref{rem.inverse}), and since $\mathcal{P}(\mathcal{P}+\lambda)^{-1}$ is bounded in $L^r(\R^d)$ uniformly in $\lambda>0$ and $\mathcal{P}(\mathcal{P}+\lambda)^{-1}h\rightarrow h$ as $\lambda\rightarrow +  0$ in $L^r(\R^d)$ for any $h\in L^r(\R^d)$, we conclude that $g_\lambda =   (\mathcal{Q} + \lambda)^{-1}M_{1/b} \nabla_x\cdot F'$ converges to $g=M_{1/b} K^*F'$ as $\lambda\rightarrow + 0$ weakly in $L^r(\R^d)$. On the other hand,  $\mathcal{Q}g_\lambda$ converges to $M_{1/b}\nabla_x\cdot F'$ strongly in $L^r(\R^d)$, which  implies  from the reflexivity of $L^r(\R^d)$ that $g\in D_{L^r}(\mathcal{Q})$ and $\mathcal{Q}g = M_{1/b} \nabla_x\cdot F'$. Hence, we have  $M_{1/b}\nabla_x\cdot F'\in {\rm Ran}(\mathcal{Q)}$ and $\mathcal{Q}^{-1}  M_{1/b} \nabla_x\cdot F' =g = M_{1/b} K^* F'$. This proves the claim.  Next we consider  $M_b \mathcal{P} \Lambda^{-1}  F_{d+1}$. As above, we take $\lambda>0$ and note that $(\Lambda +\lambda)^{-1} F_{d+1}\in D_{L^2}(\Lambda) \cap D_{L^r}(\Lambda)$ due to the assumption {\bf (i)}. Since $D_{L^2}(\mathcal{P})=D_{L^2}(\Lambda)=H^1(\R^d)$ and $M_b \mathcal{P}=\Lambda + M_{\bf{a}}\cdot \nabla_x$ by Theorem \ref{thm.factorization},  we have for any $\varphi\in C_0^\infty (\R^d)$,
\begin{align*}
\langle M_b \mathcal{P} (\Lambda + \lambda)^{-1} F_{d+1}, \varphi \rangle _{L^2(\R^d)}&  = \langle \Lambda (\Lambda + \lambda)^{-1} F_{d+1},\varphi \rangle _{L^2(\R^d)} \\
& ~~~ + \langle M_{{\bf a}}\cdot \nabla_x \Lambda^{-1}\Lambda  (\Lambda + \lambda)^{-1} F_{d+1},\varphi \rangle _{L^2(\R^d)}.
\end{align*}  
As state in Remark \ref{rem.inverse} the operator  $M_{{\bf a}}\cdot \nabla_x \Lambda^{-1}: {\rm Ran}(\Lambda)\subset L^r(\R^d) \rightarrow L^r(\R^d)$ is extended as a bounded operator $L$ in $L^r(\R^d)$. Thus we see
\begin{align*}
\| M_b \mathcal{P} (\Lambda + \lambda)^{-1} F_{d+1} \|_{L^r(\R^d)} \leq C \| \Lambda (\Lambda+\lambda)^{-1} F_{d+1}\|_{L^r(\R^d)} \leq C \| F_{d+1}\|_{L^r(\R^d)}
\end{align*}
with $C$ independent of $\lambda>0$. By the assumptions {\bf (i)} and {\bf (iii)} this estimate implies $(\Lambda + \lambda)^{-1} F_{d+1}\in D_{L^r}(\mathcal{P})$, and we have  $M_b \mathcal{P} (\Lambda + \lambda)^{-1} F_{d+1} = (I+ L )\Lambda (\Lambda+\lambda)^{-1}F_{d+1}$. Since $\Lambda (\Lambda+\lambda)^{-1}F_{d+1}\rightarrow F_{d+1}$ as $\lambda\rightarrow  + 0$ in $L^r(\R^d)$, we have $M_b \mathcal{P}\Lambda^{-1} F_{d+1} = \displaystyle \lim_{\lambda \rightarrow + 0} M_b \mathcal{P} (\Lambda+\lambda)^{-1} F_{d+1} = (I+L) F_{d+1}$ in $L^r(\R^d)$. Finally we consider $\Lambda^{-1}M_b \mathcal{Q}  F_{d+1}$. For any $\lambda, t>0$ we see 
\begin{align*}
\langle (\Lambda + \lambda)^{-1} M_b \mathcal{Q} e^{-t\mathcal{Q}}  F_{d+1}, \varphi \rangle _{L^2(\R^d)} & = \langle F_{d+1}, M_b \mathcal{P} e^{-t\mathcal{P}}  (\Lambda + \lambda)^{-1} \varphi \rangle _{L^2(\R^d)},~~~~~~~\varphi\in C_0^\infty (\R^d).
\end{align*}
As is proved above, $(\Lambda+\lambda)^{-1}\varphi \in D_{L^{r'}}(\mathcal{P})$ and we have $\| \mathcal{P}(\Lambda+\lambda)^{-1}\varphi \|_{L^{r'}(\R^d)}\leq C \|\varphi \|_{L^{r'}(\R^d)}$ with $C$ independent of $\lambda>0$ as well as $\mathcal{P}(\Lambda+\lambda)^{-1}\varphi \rightarrow M_{1/b} (I+ L) \varphi$ (as $\lambda\rightarrow + 0$) in $L^{r'}(\R^d)$. This implies that  for each $t>0$ the function $(\Lambda + \lambda)^{-1} M_b \mathcal{Q} e^{-t\mathcal{Q}}  F_{d+1}$ converges to $(I+L^*)e^{-t\mathcal{Q}} F_{d+1}$ as $\lambda\rightarrow + 0$ weakly in $L^r(\R^d)$, and from the reflexivity of $L^r(\R^d)$  we have $M_b \mathcal{Q}e^{-t\mathcal{Q}} F_{d+1}\in {\rm Ran}(\Lambda)$ in $L^r(\R^d)$ and $\Lambda^{-1} M_b \mathcal{Q}e^{-t\mathcal{Q}} F_{d+1} = (I+L^*)e^{-t\mathcal{Q}} F_{d+1}$. It is now easy to see the limit  $\displaystyle \lim_{t\rightarrow 0} \Lambda^{-1} M_b \mathcal{Q}e^{-t\mathcal{Q}} F_{d+1} = (I+L^*) F_{d+1}$ satisfies the desired estimate. The proof is complete.
\hfill
$\square$

\begin{rem}\label{rem.lem.thm.helmholtz}{\rm From the proof we have $\| \Lambda ^{-1} M_b \mathcal{Q} e^{-t\mathcal{Q}} \psi \|_{L^r(\R^d)} = \| (I + L^*) e^{-t\mathcal{Q}} \psi \|_{L^r(\R^d)} \leq C \| \psi \|_{L^r(\R^d)}$ for any $\psi \in L^r(\R^d)$ by the density argument. In particular, $M_b \mathcal{Q} e^{-t\mathcal{Q}} \psi \in {\rm Ran}(\Lambda)$ in $L^r(\R^d)$ for any $\psi \in L^r(\R^d)$. 
}
\end{rem}

\begin{lem}\label{lem.thm.helmholtz.3} Let $F=(F',F_{d+1})\in (C_0^\infty (\R^{d+1}_+))^d\times C_0^\infty (\R^{d+1}_+)$ and set $G= -(F_{d+1} + M_{{\bf a}} \cdot F')$. Let $\gamma$ be the trace operator to the boundary $\partial\R^{d+1}_+$. Assume that {\bf (i)} - {\bf (iii)} hold. Then $ \gamma G + 
M_b \int_0^\infty e^{-s \mathcal{Q}} M_{1/b} (-\nabla_x\cdot F' +\partial_s G) (s) \dd s$ belongs to ${\rm Ran} (\Lambda)$ in $L^r (\R^d)$.
\end{lem} 

\noindent {\it Proof.} By the integration by parts with respect to the time variable, we see
\begin{align}
& ~~~ \gamma G + M_b \int_0^\infty e^{-s \mathcal{Q}} M_{1/b} (-\nabla_x\cdot F' +\partial_s G) (s) \dd s \nonumber \\
& = M_b \int_0^\infty \mathcal{Q} e^{-s\mathcal{Q}} M_{1/b}  G (s) \dd s - M_b \int_0^\infty e^{-s \mathcal{Q}} M_{1/b} \nabla_x\cdot F' (s) \dd s\nonumber \\
& = M_b \mathcal{Q} \int_0^\infty e^{-s\mathcal{Q}} \big ( M_{1/b} G (s) - \mathcal{Q}^{-1} M_{1/b} \nabla_x\cdot F'(s) \big )\dd s.\label{eq.lem.thm.helmholtz.3}
\end{align}
Here we have used Lemma \ref{lem.thm.helmholtz.2}. Again from Lemma \ref{lem.thm.helmholtz.2} and Remark \ref{rem.lem.thm.helmholtz} the expression \eqref{eq.lem.thm.helmholtz.3} implies the assertion of Lemma \ref{lem.thm.helmholtz.3}. The proof is complete.
\hfill
$\square$

\begin{cor}\label{cor.lem.thm.helmholtz.3} The conclusion of Lemma \ref{lem.thm.helmholtz.3} holds for $r=2$.

\end{cor}

\noindent {\it Proof.} Proposition \ref{prop.proof.L^2} in the next section and Lemma \ref{lem.thm.helmholtz.3} prove the claim. The proof is complete.
\hfill
$\square$

\noindent {\it Proof of Proposition \ref{prop.helmholtz}.} Since $F\in (C_0^\infty (\R^{d+1}_+))^{d+1}$, there exists a unique weak solution $w\in \dot{H}^1 (\R^{d+1}_+)$ to \eqref{eq.helmholtz}. Then Corollary \ref{cor.lem.thm.helmholtz.3} and Theorem \ref{thm.representation.formula}  lead to the following representation: 
\begin{align}
 w (t)  & =  e^{-t \mathcal{P}} \Lambda^{-1} 
\big ( \gamma G + 
M_b \int_0^\infty e^{-s \mathcal{Q}} M_{1/b} 
(-\nabla_x\cdot F' +\partial_s G) (s) \dd s \big )   
\notag
\\
&~~~ + 
\int_0^t e^{-(t-s) \mathcal{P}} 
\int_s^\infty e^{-(\tau-s) \mathcal{Q}} 
M_{1/b} (-\nabla_x\cdot F' +\partial_s G) (\tau)  
\dd \tau \dd s\nonumber \\
& = e^{-t \mathcal{P}} \Lambda^{-1} M_b \mathcal{Q} \int_0^\infty e^{-s\mathcal{Q}} \big ( M_{1/b} G - \mathcal{Q}^{-1} M_{1/b} \nabla_x\cdot F' \big ) (s) \dd s\nonumber  \\
& ~ + 
\int_0^t e^{-(t-s) \mathcal{P}} \bigg ( - M_{1/b} G (s) +  \mathcal{Q}
\int_s^\infty e^{-(\tau-s) \mathcal{Q}} \big (
M_{1/b} G -\mathcal{Q}^{-1} M_{1/b} \nabla_x\cdot F') (\tau) \dd \tau  \bigg )\dd s. \label{eq.formula2}
\end{align} 
Here we have also used \eqref{eq.lem.thm.helmholtz.3}, Lemma \ref{lem.thm.helmholtz.2}, and the integration by parts. Set 
\begin{align*}
v(t) = -M_{1/b}G(t) + \mathcal{Q} \int_t^\infty e^{-(s-t) \mathcal{Q}} \big ( M_{1/b} G - \mathcal{Q}^{-1} M_{1/b} \nabla_x\cdot F' \big ) (s) \dd s. 
\end{align*}
Note that $\gamma G=0$ for $F\in (C_0^\infty (\R^{d+1}_+))^{d+1}$, and thus, we have $\gamma v = \mathcal{Q} \int_0^\infty e^{-s \mathcal{Q}} \big ( M_{1/b} G - \mathcal{Q}^{-1} M_{1/b} \nabla_x\cdot F' \big ) (s) \dd s$.  Then the solution of \eqref{eq.helmholtz} is written 
in the form 
$w=w_1 + w_2$, 
where each $w_i$ is given by 
\begin{align}
w_1 (t) = \int_0^t e^{-(t-s)\mathcal{P}} v (s) \dd s, ~~~~~~~~~~~~~~~ w_2 (t) = e^{-t\mathcal{P}}\Lambda^{-1} M_b \gamma v.
\end{align}
By the assumption \eqref{condition1-1} it suffices to estimate $\partial_t w_i$ and $\mathcal{P} w_i$.

\noindent {\it Step 1: Estimate of $w_1$.} From the maximal regularity \eqref{condition4}  and Lemma \ref{lem.helmholtz.1} we have 
\begin{align*}
\|\mathcal{P} w_1 \|_{L^q(\R_+; L^r(\R^d) )} & \leq C \| v \|_{L^q(\R_+; L^r(\R^d) )}  \\
& \leq C \big ( \| G\|_{L^q(\R_+; L^r(\R^d) )}  + \| \mathcal{Q}^{-1} M_{1/b} \nabla_x\cdot F'\|_{L^q(\R_+; L^r(\R^d) )} \big )
\end{align*}
Thus Lemma \ref{lem.thm.helmholtz.2} yields $\|\mathcal{P} w_{1} \|_{L^q (\R_+; L^r (\R^d))} \leq C_q \| F\|_{L^q (\R_+; L^r (\R^d))}$. The estimate of $\partial_t w_1$ is obtained in the same manner.

\noindent {\it Step 2: Estimate of $w_2$.}  We decompose  $\gamma v$ as $\gamma v = v_1 + v_2$, where 
\begin{align*}
v_1 & = \mathcal{Q} \int_0^t e^{-s\mathcal{Q}}\big ( M_{1/b} G - \mathcal{Q}^{-1} M_{1/b} \nabla_x\cdot F' \big ) (s)   \dd s, \\
v_2 & =  \mathcal{Q}\int_t^\infty e^{-s\mathcal{Q}} \big ( M_{1/b} G - \mathcal{Q}^{-1} M_{1/b} \nabla_x\cdot F' \big ) (s) \dd s.
\end{align*}
Motivated by this decomposition we introduce the linear operators $T_i$, $i=1,2$, defined by 
\begin{align}
T_1 [\phi] (t) & = \mathcal{P} e^{-t\mathcal{P}}  \Lambda^{-1} M_b \mathcal{Q} \int_0^t e^{-s\mathcal{Q}} \phi (s) \dd s,\\
T_2 [\phi] (t) & =   e^{-t\mathcal{P}} \mathcal{P} \Lambda^{-1}M_b  \mathcal{Q} \int_t^\infty e^{-s\mathcal{Q}} \phi (s) \dd s.
\end{align}
Each of $T_i$ makes sense for $\phi\in L^q (\R_+; L^r (\R^{d}))$ due to Lemma \ref{lem.thm.helmholtz.2} and the density argument. Clearly we have $\mathcal{P}w_{2} = \sum_{i=1,2} T_i[\phi]$ with $\phi=M_{1/b} G - \mathcal{Q}^{-1} M_{1/b} \nabla_x\cdot F' $. Lemma \ref{lem.thm.helmholtz.2}, \eqref{condition2} and \eqref{condition2'} yield  for $1<q<\infty$,
\begin{align}
\| T_1 [\phi] (t) \|_{L^r (\R^d)} \leq C t^{-1} \int_0^t \| e^{-s\mathcal{Q}} \phi (s) \|_{L^r(\R^d)} \dd s&  \leq C t^{-1} \int_0^t \| \phi(s) \|_{L^r(\R^d)} \dd s\nonumber  \\
& \leq C t^{-1/q} \| \phi \|_{L^{q} (\R_+; L^r (\R^d))},\label{proof.thm.helmholtz.1}\\
 \| T_2 [\phi] (t) \|_{L^r (\R^d)} \leq C \int_t^\infty \| \mathcal{Q}   e^{-s\mathcal{Q}}\phi (s) \|_{L^r (\R^d)} \dd s & \leq C \int_t^\infty s^{-1} \| \phi (s) \|_{L^r (\R^d)}\dd s\nonumber \\ & \leq C t^{-1/q} \| \phi \|_{L^{q} (\R_+; L^r(\R^d))}.\label{proof.thm.helmholtz.2}
\end{align}
The estimates \eqref{proof.thm.helmholtz.1}-\eqref{proof.thm.helmholtz.2} show that each $T_i$ is a bounded operator  from $L^q(\R_+; L^r(\R^d))$ to $L^{q,\infty}(\R_+; L^r(\R^d))$ for all $1<q<\infty$. Thus by the Marcinkiewicz interpolation theorem each $T_i$ is bounded in  $L^q(\R_+; L^r(\R^d))$ for all $1<q<\infty$. Hence we have arrived at $\| \mathcal{P} w_{2
}\|_{L^q(\R_+; L^r(\R^d))}\leq C \| F \|_{ L^q(\R_+; L^r(\R^d))}$ since $\|M_{1/b} G - \mathcal{Q}^{-1} M_{1/b}\nabla_x\cdot F'\|_{L^q(\R_+: L^r(\R^d))}\leq C \| F \|_{L^q(\R_+; L^r(\R^d) )}$. The estimate of $\partial_t w_2$ is the same as $\mathcal{P} w_2$. The proof of Proposition \ref{thm.helmholtz} is complete.
\hfill
$\square$

Let us  recall that the weak formulation of the Neumann problem \eqref{neumann.helmholtz} is equivalent with \eqref{eq.helmholtz.weak} through the transformation \eqref{def.pushforward}. Approximating $F\in (L^q(\R_+;L^r(\R^d))^{d+1}$ by vector fields in $(C_0^\infty (\R^{d+1}_+))^{d+1}$, we obtain from Proposition \ref{prop.helmholtz} the following 
\begin{thm}\label{thm.neumann.helmholtz}
Let $\eta$ be a globally Lipschitz function and let $\Omega$ be a domain 
given in \eqref{domain}. Assume that {\bf (i)} - {\bf (iii)} hold and let $1<q<\infty$, $\min\{m,m'\}\leq r\leq \max\{m,m'\}$. Then for any $f\in \big (Y^{q,r} (\Omega) \big )^{d+1}$ there exists  a unique (up to constant) weak solution
$p\in L^1_{loc}(\Omega)$ to the Neumann problem \eqref{neumann.helmholtz} satisfying
\begin{align}
\|\nabla p \|_{Y^{q,r}(\Omega)}\le C\|f\|_{Y^{q,r}(\Omega)}.
\end{align}
Here $C$ depends only on $d,~q,~\|\nabla_x\eta\|_{L^\infty (\R^d)}$, and the constants of the estimates in {\bf (i)} - {\bf (iii)}.  
\end{thm}

\noindent {\it Proof.} The proof of the existence proceeds as described above, and we omit the details. The uniqueness follows from the solvability of \eqref{eq.helmholtz.weak} in the adjoint space $L^{q'}(\R_+; L^{r'}(\R^d))$ and the duality. The proof is complete.
\hfill
$\square$

By the standard argument as in \cite[Lemma III.1.2]{Galdi},
we are able to show the validity of the Helmholtz decomposition in $(Y^{q,r}(\Omega))^{d+1}$ under the assumptions {\bf (i)} - {\bf (iii)}. Here we give a sketch of the proof for completeness. 
To this end we first recall a useful lemma in 
\cite[Lemma 7]{FM}, \cite[Lemma III 1.1]{Galdi}.
\begin{lem}
\label{FM}
Let $\Omega$ be a simply connected set in $\R^{d+1}$.
Suppose that $u \in (L^1_{loc}(\Omega))^{d+1}$ veryfies
$$
\int_\Omega u \cdot \varphi =0 \qquad \mathrm{for\ all} 
\quad \varphi \in C^\infty_{0,\sigma}(\Omega).
$$
Then there exists a scalar function $p \in W^{1,1}_{loc}(\Omega)$ 
such that $u=\nabla p$.
\end{lem}

\begin{thm}\label{thm.helmholtz} Under the assumptions of Theorem \ref{thm.neumann.helmholtz}, the space $(Y^{q,r}(\Omega) )^{d+1}$ admits the Helmholtz decomposition. 
\end{thm}

\noindent {\it Proof.}\,  For $f \in (Y^{q,r}(\Omega))^{d+1}$ let $p\in L^1_{loc} (\Omega)$ be the weak 
solution to  \eqref{neumann.helmholtz}  given in 
Theorem \ref{thm.neumann.helmholtz}  which satisfies $\nabla p \in Y^{q,r}_G(\Omega)$, and 
set $u=f-\nabla p.$  Then, we have $u \in (Y^{q,r}_{G}(\Omega))^{\perp}$
where $X^\perp$ denotes the annihilator of the set $X$.
On the other hand, by Lemma \ref{FM}, we see
$$
(Y^{q,r}_{\sigma})^{\perp} \subset Y^{q',r'}_G, \ 
\mathrm{or \ equivalently}, \ Y^{q,r}_{\sigma} \supset (Y^{q',r'}_G)^{\perp}
$$
Thus we have proved that $u \in Y^{q,r}_{\sigma}$. 

It remains to show that the representation $f=u+\nabla p$ 
is unique. This is equivalent to show that the equality
\begin{equation}
u=\nabla p, \quad u \in Y_\sigma^{q,r}, \quad \nabla p \in Y_G^{q,r}
\end{equation}
holds if and only if $u\equiv \nabla p \equiv 0$. 
To this end, we observe that 
$$
Y_\sigma^{q,r} \subset (Y_G^{q',r'})^\perp,
$$
which implies
$$
\ \qquad \int_\Omega \nabla p \cdot \nabla \varphi=0 
\qquad \mathrm{for\ all} 
\quad \varphi \in 
\{f\in L^1_{loc}(\Omega)\ |\ \nabla f \in (Y^{q',r'}(\Omega))^{d+1} \}.
$$
By the uniqueness of weak solutions to \eqref{neumann.helmholtz},
we have $u \equiv \nabla p \equiv 0$. 
The proof is complete.
\hfill
$\square$

\subsection{Proof of Theorem \ref{thm.helmholtz.intro}}\label{subsec.helmholtz.L^2}

In this section we prove the next proposition, which immediately leads to Theorem \ref{thm.helmholtz.intro} thanks to Theorem \ref{thm.helmholtz} in the previous section.
\begin{prop}\label{prop.proof.L^2} When $m=m'=r=2$, the assumptions {\bf (i)} - {\bf (iii)} in Section \ref{subsec.helmholtz.L^r} are valid.
\end{prop} 

\noindent {\it Proof.}  It is well-known that $\Lambda$ is self-adjoint in $L^2(\R^d)$, and \eqref{def.Lambda} implies the boundedness of $\{e^{-t\Lambda}\}_{t\geq 0}$ in $L^2 (\R^d)$. We also observe from Theorem \ref{thm.factorization} that $D_{L^2}(\mathcal{P}) = D_{L^2}(\Lambda) = H^1 (\R^d)$ and  $\mathcal{P}$ generates a strongly continuous and bounded analytic semigroup in $L^2(\R^d)$. Then it suffices to check the estimates \eqref{condition1-1}, \eqref{condition1-2}, and \eqref{condition3}  for $r=2$. In fact, these estimates are already known as a consequence of the Rellich type identity, which is a classical tool in the  study of the Dirichlet and Neumann problems for the Laplace equations; see Remark \ref{rem.rellich} below for references. Here we give the detailed proof of them including \eqref{condition2'} in order to make this paper self-contained as much as possible.

\noindent {\it Step 1: Proof of {\rm {\bf (ii)}}.} We will prove 
\begin{align}
& \| M_{\sqrt{b}} \mathcal{P} \varphi \|_{L^2(\R^d)}  =  \| \nabla_x \varphi \|_{L^2(\R^d)} \label{proof.thm.helmholtz.2.3}\\
&C_1 \| \nabla_x \varphi \|_{L^2(\R^d)}\leq  \|  \Lambda \varphi \|_{L^2(\R^d)} \leq C_2 \|\nabla_x \varphi \|_{L^2(\R^d)},  \label{proof.thm.helmholtz.2.1}
\end{align}
where $C_1$ and  $C_2$ are positive constants  depending only on $\|\nabla_x\eta \|_{L^\infty (\R^d)}$. Let us recall that in our case Theorem \ref{thm.factorization}  holds with $\mathcal{A}'=-\Delta_x$. Thus   \eqref{eq.factorization'} yields \eqref{proof.thm.helmholtz.2.3}. Next by the relation $\mathcal{P}=M_{1/b} \Lambda - M_{(\nabla_x\eta)/b}\cdot \nabla _x$ the right-hand side of \eqref{proof.thm.helmholtz.2.3} is written as 
\begin{align}
\|M_{\sqrt{b}} \mathcal{P} \varphi \|_{L^2(\R^d)}^2 
& = 
\| M_{1/\sqrt{b}}\Lambda \varphi - M_{(\nabla_x\eta)/\sqrt{b}} \cdot \nabla_x \varphi  \|_{L^2(\R^d)}^2 \nonumber \\
& = \|M_{1/\sqrt{b}} \Lambda \varphi \|_{L^2(\R^d)}^2 -2\langle \Lambda \varphi, M_{(\nabla_x\eta)/b}\cdot \nabla_x \varphi  \rangle_{L^2(\R^d)} + \| M_{(\nabla_x\eta/\sqrt{b})} \cdot \nabla_x \varphi \|_{L^2(\R^d)}^2.\label{proof.thm.helmholtz.2.4}
\end{align}
Thus \eqref{proof.thm.helmholtz.2.3} and  \eqref{proof.thm.helmholtz.2.4} immediately yield
\begin{equation}
\| M_{1/\sqrt{b}} \Lambda \varphi \|_{L^2(\R^d)}^2 \leq 2 \|\nabla_x \varphi \|_{L^2(\R^d)}^2.\label{proof.thm.helmholtz.2.5}
\end{equation}
While, we derive from \eqref{proof.thm.helmholtz.2.4} that
\begin{align*}
\|M_{\sqrt{b}} \mathcal{P} \varphi  \|_{L^2(\R^d)}^2 & \leq (1+\epsilon^{-1}) \|M_{1/\sqrt{b}} \Lambda \varphi \|_{L^2(\R^d)}^2 + (1+\epsilon ) \int_{\R^d} \frac{|\nabla_x\eta |^2}{1+|\nabla_x\eta|^2} |\nabla_x \varphi |^2 \dd x,~~~~~~~~\epsilon>0.
\end{align*}
Hence, combining this with \eqref{proof.thm.helmholtz.2.3} implies 
\begin{align*}
(1+\epsilon ^{-1})\|M_{1/\sqrt{b}} \Lambda \varphi \|_{L^2(\R^d)}^2 \geq \int_{\R^d} \frac{1-\epsilon |\nabla_x\eta |^2}{1+|\nabla_x\eta|^2} |\nabla_x \varphi |^2 \dd x\geq c \|\nabla_x \varphi \|_{L^2(\R^d)}^2,
\end{align*}
if $\epsilon>0$ is small enough. Here $c>0$ depends only on $\|\nabla_x\eta\|_{L^\infty(\R^d)}$. Thus {\bf (ii)} is proved.

\noindent {\it Step 2: Proof of {\rm {\bf (iii)}}.} From Remark \ref{rem.thm.helmholtz.1.2} it suffices to show \eqref{condition3'} with $r=2$.  Set $w(t) = \int_0^t e^{-(t-s)\mathcal{P}} \phi (s) \dd s$. Then $w$ is the solution to $\partial _t w + \mathcal{P} w =\phi$ in $(t,x)\in \R_+ \times \R^d$ and $\gamma w=0$ on $\partial\R^{d+1}_+$. Hence we have from \cite[Theorem 1.10, Remark 1.11]{MM},
\begin{align}
\langle A \nabla w, \nabla w\rangle _{L^2(\R^{d+1}_+)} & =  \langle \phi, M_b (\partial_t + \mathcal{P} )w\rangle _{L^2(\R^{d+1}_+)}.\label{proof.thm.helmholtz.2.2}
\end{align}
Here we have used the boundary condition $w=0$ on $t=0$. By the uniform ellipticity the left-hand side of \eqref{proof.thm.helmholtz.2.2} is bounded from below by $c\|\nabla w \|_{L^2(\R^{d+1}_+)}^2$ with some $c>0$ depending only on $d$ and $\|\nabla_x\eta \|_{L^\infty (\R^d)}$.  On the other hand,  the 
right-hand side of \eqref{proof.thm.helmholtz.2.2} is calculated as 
\begin{align*}
{\rm R.H.S.~of~} \eqref{proof.thm.helmholtz.2.2} & =  \langle M_b \phi, \phi \rangle_{L^2 (\R^{d+1}_+)}  \leq C\| \phi \|_{L^2(\R^{d+1}_+)}^2. 
\end{align*}
This complete the proof of {\bf (iii)}.

Finally, we also give the proof of \eqref{condition2'} for $\{e^{-t\mathcal{P}}\}_{t\geq 0}$ for reader's convenience. It suffices to consider the case  $\varphi \in C_0^\infty(\R^d)$. Set $u_0 (t) = e^{-t\mathcal{P}} \varphi$. Then  we see from \eqref{condition3} with $r=2$ (which is proved in Step $2$ above) that 
\begin{align}
\| u_0 \|_{L^2(0,T; L^2 (\R^d))} \leq C T^{1/2} \| \varphi \|_{L^2(\R^d)},~~~~~~~~~~~~ T>0,\label{proof.step3.1}
\end{align}
where $C>0$ is independent of $T$ and $\varphi$.  Set $u_k (t) = (t\mathcal{P})^k u_0 (t)$ for $k=1,$ $2$. Then $u_k$ satisfies $\partial_t u_k - \mathcal{P} u_k = k \mathcal{P} u_{k-1}$ with the zero initial data.  Hence we have the representation 
$u_k(t) = k \int_0^t e^{(t-s)\mathcal{P}} \mathcal{P} u_{k-1} \dd s = k \mathcal{P} \int_0^t e^{(t-s)\mathcal{P}} u_{k-1} \dd s$. 
Thus \eqref{condition3} with $r=2$ yields 
\begin{align}
\| u_k \| _{L^2(0,T; L^2 (\R^d))} \leq C \| u_{k-1} \|_{L^2 (0,T; L^2 (\R^d ))}\leq C T^\frac12 \| \varphi\|_{L^2(\R^d)},~~~~~~~~~~T>0. 
\label{proof.step3.2}
\end{align}
In particular, \eqref{proof.step3.1} and \eqref{proof.step3.2} with $k=1$ imply that for any $T>0$ there exists $T_1\in [T/2,T]$ such that 
$\|u_1 (T_1) \|_{L^2(\R^d)} \leq C \| \varphi \|_{L^2(\R^d)}$ with $C$ independent of $T$ and $\varphi$. 
Next we calculate the evolution of  $\|M_{\sqrt{b}} u_1 (t) \|_{L^2(\R^d)}^2$ to get 
\begin{align*}
\frac{\dd}{\dd t} \| M_{\sqrt{b}} u_1 (t) \|_{L^2(\R^d)} ^2 
& = 
- 2 \langle u_1,  \Lambda u_1 - M_{\nabla_x\eta}\cdot \nabla_x u_1\rangle _{L^2(\R^d)} 
+ 2 
\langle u_1, M_b\mathcal{P} u_0 \rangle _{L^2(\R^d)}
\\
& \leq 2\langle u_1, M_{\nabla_x\eta} \cdot \nabla_x u_1\rangle _{L^2(\R^d)} + 2 t^{-1} \| M_{\sqrt{b}} u_1 (t) \|_{L^2(\R^d)}^2
\\
& \leq C \| u_1 (t) \|_{L^2(\R^d)} \| \mathcal{P} u_1 (t) \|_{L^2(\R^d)} + 2 t^{-1} \| M_{\sqrt{b}} u_1 (t) \|_{L^2(\R^d)}^2\\
&\leq C t \| \mathcal{P} u_1 (t) \|_{L^2(\R^d)}^2 + C t^{-1} \| M_{\sqrt{b}} u_1 (t) \|_{L^2(\R^d)}^2.
\end{align*}
Here we have used the definition of $u_k$, $\langle u_1, \Lambda u_1 \rangle _{L^2(\R^d)} \geq 0$, and \eqref{condition1-1} with $r=2$. Then by integrating over $[T_1,T]$ and by using the Gronwall inequality, we arrive at 
\begin{align*}
\| M_{\sqrt{b}} u_1 (T) \|_{L^2(\R^d)} ^2 & \leq C \|  M_{\sqrt{b}} u_1 (T_1) \|_{L^2(\R^d)} ^2 + C\int_{T_1}^T t \| \mathcal{P} u_1 (t) \|_{L^2(\R^d)}^2 \dd t\\
& \leq  C \|  u_1 (T_1) \|_{L^2(\R^d)} ^2 + CT_1^{-1}\| u_2  \|_{L^2 (T_1,T; L^2(\R^d ))}^2 \leq C \| \varphi \|_{L^2(\R^d)}^2.
\end{align*}
Since $T>0$ is arbitrary, we have proved
\begin{align}
\| \mathcal{P} e^{-t\mathcal{P}} \varphi \|_{L^2 (\R^d)}\leq C t^{-1} \| \varphi \|_{L^2 (\R^d)},~~~~~~~t>0.\label{proof.step3.3}
\end{align} 
Now from \eqref{proof.step3.1} there is $T_2\in [T/2,T]$ such that $\| u_0(T_2) \|_{L^2 (\R^d)}\leq C \| \varphi \|_{L^2 (\R^d)}$. Combining this with \eqref{proof.step3.3} and the equality $u_0(T)= u_0(T_2) -\int_{T_2}^T \mathcal{P} e^{-s\mathcal{P}} \varphi \dd s$, we obtain \eqref{condition2'}.  The proof of Proposition \ref{prop.proof.L^2} is complete. 
\hfill
$\square$

\begin{rem}\label{rem.rellich}
{\rm As mentioned above, for the case $r=2$, 
the coercive estimates for $\mathcal{P}$ and $\Lambda$ as in \eqref{condition1-1} - \eqref{condition1-2}
are obtained from a variant of the Rellich identity \cite{Rellich}. These estimates are used to study the behavior of harmonic functions near the boundary \cite{PW,JK, JK2}. In particular, it works even for nonsmooth domains, and in \cite{JK,JK2} the Rellich type identity was used in solving  the Dirichlet and Neumann problems in bounded Lipschitz domains. We also note that  \eqref{condition1-1} and \eqref{condition1-2} are the key to obtain  the characterization $D_{L^2} (\mathcal{P})=H^1 (\R^d)$ and $D_{L^2}(\Lambda)=H^1 (\R^d)$ for  real symmetric  but nonsmooth $A$. For a matrix $A$ of the form in Section \ref{sec.proof}, called the Jacobian type, the relation $D_{L^2}(\mathcal{P})=H^1 (\R^d)$ is proved in \cite{Dahlberg2}. The relation $D_{L^2}(\Lambda)=H^1 (\R^d)$ is related with the solvability of the Neumann problem for $L^2$ boundary data, and it is solved by \cite{JK}  in bounded Lipschitz domains. For results in more general class of $A$ including real symmetric or Hermite ones, see \cite{JK2,KKPT,AAH,AAM,AAAHK,MM} and references therein.

}
\end{rem}

\medskip

\subsection{Concluding remark}
Recent works \cite{AT,GHHS} revealed that the solvability
of the weak Neumann problem in $L^p$  ensures the analyticity of the Stokes
semigroup in $L^p$ even if the boundary in noncompact. In view of their results, it is expected that
the Stokes semigroup is analytic in the space $Y^{q,2}$ at least when the boundary is smooth enough. 
We will address this question in the forthcoming work.

\appendix

\section{Appendix}

\subsection{Semigroup $\{e^{-t\Lambda}\}_{t\geq 0}$ in $L^r(\R^d)$ for $r\in (1,\infty)$}

\begin{prop}\label{prop.DN.L^r} Let $r\in (1,\infty)$. Then the restrictions of  $\{e^{-t\Lambda }\}_{t\geq 0}$ on $L^2(\R^d)\cap L^r(\R^d)$ is extended as a  strongly continuous and bounded semigroup in $L^r (\R^d)$.
\end{prop}

\noindent {\it Proof.} Here we give only a sketch of the proof. We first consider the case $r\in [2,\infty)$. Set $u(t) = e^{-t\Lambda} f$, $f\in H^1 (\R^d) \cap L^r (\R^d)$, and set $v(s;t) = e^{-s\mathcal{P}} u(t)$, $s\geq 0$. Then we have 
\begin{align*}
\frac{\dd}{\dd t} \| u(t) \|_{L^r(\R^d)}^r & = - r \langle \Lambda u(t), |u(t)|^{r-2} u(t) \rangle _{L^2(\R^d)} \\
& = - r \int_0^\infty \langle A\nabla v(s;t), \nabla \big ( | v(s;t)|^{r-2} v(s;t) \big ) \rangle _{L^2 (\R^d)} \dd s \\
& \leq -c r  \| \nabla ( |v(\cdot; t) |^{\frac{r}{2}} ) \|_{L^2 (\R^{d+1}_+)}^2\leq -c r \| | u(t) |^{\frac{r}{2}} \|_{\dot{H}^\frac12(\R^d)}^2 .
\end{align*}
In particular, we have $\| e^{-t\Lambda} f \|_{L^r(\R^d)}\leq \| f \|_{L^r(\R^d)}$ when $r\in [2,\infty)$, and thus, when $r\in [2,\infty]$ (by taking the limit $r\rightarrow \infty$). By the dual relation $\langle e^{-t\Lambda} f, g\rangle_{L^2 (\R^d)} = \langle f, e^{-t\Lambda}g \rangle _{L^2 (\R^d)}$ we have this uniform bound also for $r\in [1,2]$. Hence, by the density argument $\{e^{-t\Lambda}\}_{t\geq 0}$ is extended as a bounded semigroup acting on $L^r(\R^d)$ for all $r\in [1,\infty]$ (note that this uniform bound holds also for $r=1,\infty$). As for the strong continuity, let $r\in (2,\infty)$, and for any $f\in L^r(\R^d)$ we take $\{f_n\}\subset C_0^\infty (\R^d)$ such that $f_n\rightarrow f$ as $n\rightarrow \infty$ in $L^r(\R^d)$. Then we see
\begin{align*}
\| e^{-t\Lambda} f - f\|_{L^r (\R^d)} &\leq \| e^{-t\Lambda} (f-f_n)\|_{L^r(\R^d)} + \| f-f_n \|_{L^r(\R^d)} + \| e^{-t\Lambda} f_n-f_n\|_{L^r(\R^d)}\\
& \leq 2 \| f -f_n \|_{L^r(\R^d)} + \|  e^{-t\Lambda} f_n-f_n\|_{L^2(\R^d)}^{\frac{2}{r}} \| e^{-t\Lambda} f_n-f_n\|_{L^\infty (\R^d)}^{1-\frac{2}{r}}\\
& \leq 2 \| f -f_n \|_{L^r(\R^d)} + 2^{1-\frac{2}{r}} \|  e^{-t\Lambda} f_n-f_n\|_{L^2(\R^d)}^{\frac{2}{r}} \| f_n \|_{L^\infty (\R^d)}^{1-\frac{2}{r}}.
\end{align*}
Since we have already known that $\{e^{-t\Lambda}\}_{t\geq 0}$ is strongly continuous in $L^2 (\R^d)$, the last estimate implies the strong continuity in $L^r(\R^d)$ for $r\in (2,\infty)$. The case $r\in (1,2)$ is proved in the same manner. The proof is complete.
\hfill
$\square$

\subsection{Semigroup $\{e^{-t\mathcal{P}}\}_{t\geq 0}$ in $L^r(\R^d)$ for $r\in [2,\infty)$}

\begin{prop}\label{prop.Poisson.L^r} Let $r\in [2,\infty)$. Then the restrictions of  $\{e^{-t\mathcal{P} }\}_{t\geq 0}$ on $L^2(\R^d)\cap L^r(\R^d)$ is extended as a  strongly continuous and bounded semigroup in $L^r (\R^d)$.
\end{prop}

\noindent {\it Proof.} Again we give only a sketch of the proof. Let $f\in L^2 (\R^d)\cap L^\infty (\R^d)$ and set $u(t) = e^{-t\mathcal{P}} f$. Then, since $u$ satisfies $\mathcal{A}u=0$ in $\R^{d+1}_+$, the maximum principle implies that 
\begin{align*}
\| u(t) \|_{L^\infty(\R^d)} \leq \| u \|_{L^\infty (\R^{d+1}_+)} \leq \| f \|_{L^\infty (\R^d)}.
\end{align*}
This estimate gives the boundedness of $e^{-t\mathcal{P}}$ in $L^\infty(\R^d)$. Since $e^{-t\mathcal{P}}$ is bounded in $L^2 (\R^d)$, the interpolation inequality yields the boundedness of $e^{-t\mathcal{P}}$ in $L^r(\R^d)$ for each $r\in (2,\infty)$. The strong continuity in $L^r(\R^d)$ is shown as in the proof  of Proposition \ref{prop.DN.L^r}. The proof is complete.  
\hfill
$\square$

\end{document}